\documentclass{fundam}

\usepackage{latexsym}
\usepackage{amssymb}
\usepackage{amsmath}
\usepackage{graphicx}
\usepackage[all]{xy}
\usepackage[T1]{fontenc}
\usepackage[english]{babel}
\usepackage{caption} 
\usepackage{mathabx} 
\usepackage{tikz}
\usepackage{caption, subcaption}
\usepackage[misc]{ifsym}

\def\I{\mathop{\mbox{\rm I}}}
\def\KK{\mathbb{K}}
\def\Ext{\mathop{\mbox{\rm Ext}}}
\def\Int{\mathop{\mbox{\rm Int}}}
\def\neu#1{\textbf{#1}}
\def\ddh{\downdownharpoons}
\def\uuh{\upupharpoons}

\begin{document}

\setcounter{page}{217}
\publyear{22}
\papernumber{2147}
\volume{188}
\issue{4}

\finalVersionForARXIV

 \title{Formal Concepts and Residuation on Multilattices}

\author{Blaise B. Koguep Njionou\\
Department of  Mathematics and Computer Science\\
University of Dschang, BP 67, Cameroon \\
blaise.koguep@univ-dschang.org
\and
Leonard Kwuida\thanks{Address for correspondence:
                    Bern University of Applied Science, Bern, Switzerland. \newline \newline
                    \vspace*{-6mm}{\scriptsize{Received January 2022; \ accepted  April  2023.}}}
\\
Bern University of Applied Sciences, Bern, Switzerland\\
leonard.kwuida@bfh.ch
\and
Celestin Lele\\
Department of  Mathematics and Computer Science\\
University of Dschang, BP 67,  Cameroon\\
celestinlele@yahoo.com}

 \maketitle

\runninghead{B.B. Koguep Njionou et al.}{Formal Concepts and Residuation on Multilattices}

\begin{abstract}
Multilattices are generalisations of lattices introduced by Mihail Benado in \cite{Bm55}. He replaced the existence of unique lower (resp. upper) bound by the existence of maximal lower (resp. minimal upper) bound(s). A multilattice will be called pure if it is not a lattice. Multilattices could be endowed with a residuation, and therefore used as set of truth-values to evaluate elements in fuzzy setting. In this paper 
we exhibit the smallest pure multilattice and show that it is a sub-multilattice of any pure multilattice.  We also prove that any bounded residuated multilattice that is not a residuated lattice has at least seven elements. We apply the ordinal sum construction to get more examples of residuated multilattices that are not residuated lattices. We then use these residuated multilattices to evaluate objects and attributes in formal concept analysis setting, and describe the structure of the set of corresponding formal concepts. More precisely, if $\mathcal{A}_i: =(A_i,\le_i,\top_i,\odot_i,\to_i,\bot_i)$, $i=1,2$ are two  complete residuated multilattices, $G$ 
and $M$ 
two nonempty sets and $(\varphi, \psi)$ a Galois connection between $A_1^G$ and $A_2^M$ that is compatible with the residuation, then we show that
 \[\mathcal{C}: =\{(h,f)\in A_1^G\times A_2^M; \varphi(h)=f \text{ and } \psi(f)=h \}\]
can be endowed with a complete residuated multilattice structure. This is a generalization of a result by Ruiz-Calvi{\~n}o and Medina
 \cite{RM12} saying that if the (reduct of the) algebras $\mathcal{A}_i$, $i=1,2$ are complete multilattices,  then $\mathcal{C}$ is a complete multilattice.
 \end{abstract}


\section{Introduction}
In the theory of fuzzy concept analysis, the underlying set of truth values is generally a lattice. Medina and Ruiz-Calvi{\~n}o proposed in \cite{RM12} a new approach of fuzzy concept analysis by using multilattices,  introduced by Benado~\cite{Bm55}, as underlying set of truth values. The main idea is to relax the requirement on the existence of least upper bounds and greatest lower bounds, and ask only that the set of minimal upper bounds and the set of maximal lower bounds of any pair of elements are non-empty. In this paper, we show that residuation on multilattices can be useful in fuzzy concept analysis to evaluate attributes and objects. Especially, we will use the adjoint pair in a residuated multilattice to build a residuated concept  multilattice.

We organize the present contribution as follows: Section~\ref{s:CLmL} recalls some notions in Formal Concept Analysis (FCA) and multilattices. Section~\ref{s:TDS} 
shows how to build some residuated multilattices. The idea here is to prove that there are enough residuated multilattices, that can be latter used for evaluating objects and attributes. We start by 
proving that every bounded pure residuated multilattice has at least seven elements. Inspired by what was done on lattices (\cite{EMM13,Em20}), we apply the ordinary sum construction to create new residuated multilattices from  old ones. We briefly recall how  Medina and Ruiz-Calvi{\~n}o used multilattices as truth degree sets. In Section~\ref{s:RCML}, we use residuated multilattices as the set of truth-values to evaluate attributes and objects, and find a condition under which the set of all concepts forms a residuated multilattice.

\section{Concept lattices and multilattices} \label{s:CLmL}

A \neu{formal context} is set up with the sets $G$ of objects, $M$  of attributes and a binary relation
$\I\subseteq G\times M$. We denote it  by $\KK:=(G,M,\I)$, and write $(g,m)\in\I$ or $g\I m$ to mean that the
object $g$ has the attribute $m$. In the whole paper, we assume that $G$ and $M$ are non-empty. To extract knowledge from formal contexts, one can get clusters of objects/attributes  called concepts.
The Port-Royal Logic School considers a \neu{concept} as defined by two parts: an extent and an intent.
The \neu{extent} contains all entities belonging to the concept and the \neu{intent} is the set of all attributes common to all entities in the concept. To formalize the notion of concept, we need the
\neu{derivation operator} $'$, defined on $A\subseteq G$ and $B\subseteq M$ by: 
	\begin{align*}
	A':=\{m\in M; g\I m\text{ for all }g\in A\} \quad \text{ and } \quad
	B':=\{g\in G; g\I m\text{ for all }m\in B\}.
	\end{align*}
A \neu{formal concept} is a pair $(A,B)$ with $A'=B$ and $B'=A$. The set of formal concepts of $\KK$ is denoted by $\mathfrak{B}(\KK)$. It forms a complete lattice, called \neu{concept lattice} of the context $\KK$, when ordered by the \neu{concept hierarchy} below:
$$\left(A_{1},B_{1}\right)\leqslant\left(A_{2},B_{2}\right):\iff A_{1}\subseteq A_2.$$

Recall that a \neu{lattice} is a poset in which all finite subsets have a least upper bound and a greatest lower bound. A lattice is complete if every subset has a least upper bound and a greatest lower bound.
For any subset $X$ of $L$, we denote by $\bigvee X$ its least upper bound and by $\bigwedge X$ its greatest lower bound, whenever they exist. For $X=\{x,y\}$ we write $x\vee y:=\bigvee X$ and $x\wedge y:=\bigwedge X$.

\medskip
For any context $(G,M,I)$ the pair $(',')$ forms a Galois connection between $\mathcal{P}(G)$ and $\mathcal{P}(M)$ and $c:X\mapsto X''$ a closure operator (on $\mathcal{P}(G)$ or $\mathcal{P}(M)$, ordered by the inclusion). Recall that a \neu{closure operator} on a poset $(P,\le)$ is a map $c:P\to P$ such that
\[x\le c(y) \iff c(x)\le c(y)\quad \text{ for all } x,y\in P.\]
A pair $(\varphi,\psi)$ is a \neu{Galois connection} between the two posets $(P,\le)$ and $(Q,\le)$ if
$\varphi: P\to Q$ and $\psi: Q\to P$ are maps such that
$$p\le \psi q \iff q\le \varphi p\quad \text{ for all } p\in P \text{ and  for all } q\in Q.$$
In general we will call any pair $(p,q)$ a \neu{concept} if $\varphi p=q$ and $\psi q=p$, whenever $(\varphi,\psi)$ is a Galois connection.
The maps $\varphi\circ\psi$ and $\psi\circ \varphi$ (denoted by $\varphi\psi$ and $\psi\varphi$ for short) are the corresponding closure operators on $P$ and $Q$ respectively. An element $p\in P$ is \neu{closed} if $\psi\varphi(p)=p$. Similarly $q\in Q$ is closed if $\varphi\psi(q)=q$.

\medskip
Objects, attributes or incidences can be of fuzzy nature. Several sets of truth degrees have been considered,
starting from the interval $[0,1]$ \cite{ZL65} to lattices \cite{GJ67}, multilattices \cite{MOR13} and residuations on these structures.
%
 In \cite{RM12} Ruiz-Calvi{\~n}o and Medina investigate the use of multilattices as underlying set of
 truth-values for attributes and objects, and prove that the set of all concepts is a multilattice.
 We will use  residuated multilattices as the set of truth-values to evaluate attributes and objects,
 and show that the set of all concepts in this case is a residuated multilattice. 

\medskip
 Multilattices generalize lattices and were introduced by Mihail Benado~\cite{Bm55}. Let $(P,\le)$ be a
 poset and $x,y\in P$. We say that $x$ \neu{is below} $y$ or  $y$ \neu{is above} $x$, whenever
 $x\le y$.\footnote{We also accept that any element is below itself.}

\begin{definition}\cite{RM12}
 A poset $(P,\le)$ is called \neu{coherent} if every chain has supremum and infimum.
\end{definition}

 \begin{definition}\cite{MOR07}
 A poset $(P,\le)$ is called a \neu{multilattice} 
 if for any finite subset $X$ of $P$,
 each upper bound of $X$ is above a minimal upper bound of $X$ and each lower bound of $X$ is below a maximal
 lower bound of $X$. 	
\end{definition}
For any subset $X$ of $P$, we denote by $\sqcup X$  (resp. $\sqcap X$) the set of its minimal upper bounds
(resp. maximal lower bounds). In general, $\sqcup X$ and $\sqcap X$ could be empty.  

If $(P,\le)$ is a bounded finite poset, then $\sqcup X$ and $\sqcap X$ are non-empty. For $X=\{x,y\}$ we write $x\sqcup y$  and $x\sqcap y$ instead of $\sqcup X$  and $\sqcap X$. When  $\sqcup X$  or $\sqcap X$ is a singleton it will be considered as an element of $P$, i.e., we write $\sqcup X\in P$  or $\sqcap X\in P$ instead of $\sqcup X\subseteq P$  or $\sqcap X\subseteq P$. If $(P,\le)$ is a lattice and $X$ a finite subset of $P$, then the set of upper (resp. lower) bounds of $X$ is non-empty and has exactly one minimal (resp. maximal) element, namely $\vee X$ (resp. $\wedge X$). Moreover, any upper (resp. lower) bound of $X$ is above (resp. below) $\vee X$ (resp. $\wedge X$). Thus any lattice is a multilattice with $\sqcup X=\{\vee X\}$ and $\sqcap X=\{\wedge X\}$. We call a multilattice \neu{pure} if it is not a lattice, and \neu{full} if $x\sqcap y$ and $x\sqcup y$ are non-empty for all $x,y$. 

 FCA can be seen as applied theory of complete lattices. The completeness can be carried out to multilattices as follows:

\begin{definition}\cite{MOR07}
A poset $(P, \le)$ is a \neu{complete multilattice}~\footnote{In \cite{RM12} the authors defined
complete multilattices as coherent posets $(P,\le)$ with no infinite antichain and
$\sqcup X\neq \emptyset\neq\sqcap X$ for all $X\subseteq P$.} if for all $X\subseteq P$, the
 sets $\sqcup X, \sqcap X$ are non-empty and each upper bound of $X$ is above an element of $\sqcup X$
 and each lower bound is below an element of $\sqcap X$.
\end{definition}

\begin{figure}[!htbp]
\vspace*{-3mm}
	\begin{center}
	\begin{tikzpicture}[scale=0.7]
	\draw[fill] (2,0) circle (0.05);  
	\draw[fill] (1,1) circle (0.05);  
	\draw[fill] (3,1) circle (0.05);  
	\draw[fill] (3,3) circle (0.05);  
	\draw[fill] (1,3) circle (0.05);  
	\draw[fill] (2,4) circle (0.05);  
	\draw (2,0) -- (1,1) -- (1,3) -- (2,4)--(3,3)--(3,1)--(2,0);
	\draw (1,1) -- (3,3);
	\draw (1,3) -- (3,1);
	\draw (2,0) node[below]{$\bot$};  
	\draw (1,1) node[left]{a};  
	\draw (3,1) node[right]{b};  
	\draw (3,3) node[right]{d};  
	\draw (1,3) node[left]{c};  
	\draw (2,4) node[above]{$\top$};  
	\end{tikzpicture}
\hspace{1cm}
	\begin{tikzpicture}[scale=0.7]
	\draw[fill] (2,0) circle (0.05);  
	\draw[fill] (1,1) circle (0.05);  
	\draw[fill] (3,1) circle (0.05);  
	\draw[fill] (3,3) circle (0.05);  
	\draw[fill] (1,3) circle (0.05);  
	\draw (3,1)--(2,0)--(1,1) -- (1,3) -- (3,1) -- (3,3) -- (1,1);
	\draw (2,0) node[below]{$\bot$};  
	\draw (1,1) node[left]{a};  
	\draw (3,1) node[right]{b};  
	\draw (3,3) node[right]{d};  
	\draw (1,3) node[left]{c};  
	\end{tikzpicture}
\end{center}\vspace*{-5mm}
\caption{A complete and pure multilattice (left) and a non-complete multilattice (right).} \label{fig:cmlat}\vspace*{-2mm}
\end{figure}
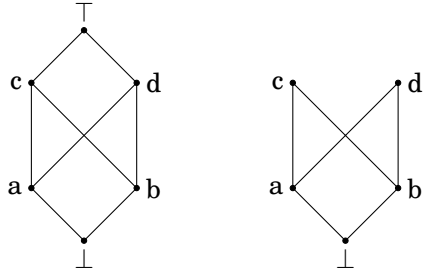

On Figure~\ref{fig:cmlat} we have the smallest bounded poset that is a complete multilattice but not a lattice. We will usually refer to it as the multilattice $\mathrm{ML_6}$. All finite lattices are complete lattices, and all finite bounded posets are complete multilattices.

In the framework of multilattices, the concept of homomorphism has been originally introduced by M.
Benado \cite{Bm55}.

\begin{definition}\cite{Bm55}
A map $h: M \rightarrow P$ between two multilattices $M$ and $P$ is said to be a \textbf{homomorphism} if
$h(x\sqcup y)\subseteq h(x)\sqcup h(y)$ and $h(x\sqcap y) \subseteq h(x) \sqcap h(y)$, for all $x, y\in M$.
\end{definition}
When the initial multilattice is full, the notion of homomorphism can be characterized in terms of equalities.

\begin{proposition}\cite{CCMO14}
Let $h: M \rightarrow P$ be a map between multilattices where $M$ is full. Then $h$ is a homomorphism if and only if, for all $x$, $y\in M$,
 $h(x\sqcup y)= (h(x)\sqcup h(y))\cap h(M)$ and $h(x\sqcap y) = (h(x) \sqcap h(y))\cap h(M)$.
\end{proposition}
A homomorphism $h$ will be called an \textbf{isomorphism} when it is a bijection.

\medskip
In \cite{MOR07}, J.Medina, M. Ojeda-Aciego and J. Ruiz-Calvi{\~n}o defined two different types of submultilattices, namely : full submultilattice (or f-submultilattice) and restricted submultilattice
(or r-submultilattice).

\begin{definition}\cite{MOR07}
Let $(P,\le)$ be a multilattice and $X$ be a nonempty subset of $P$.
\begin{itemize}
\item[(i)] $X$ is called a \textbf{full submultilattice} ($f$-\textbf{submultilattice}) of $P$ if for all $x, y\in X$, $x\sqcup y\subseteq X$ and $x\sqcap y\subseteq X$ (SML-1).
     \item[(ii)] $X$ is called a \textbf{restricted submultilattice} (\textbf{$r$-submultilattice}) of $P$ if for all $x, y\in X$, $(x\sqcup y)\cap X\ne \emptyset$ and $(x\sqcap y)\cap X\ne \emptyset$ (SML-2).
\end{itemize}
\end{definition}
It was proved in \cite{MOR07} that if $X$ is a full submultilattice of a multilattice $P$, then equipped with
the restriction of the partial order from  $P$, $X$ is a multilattice on its own right. However, this is not
the case for restricted submultilattices as we can see in the following example.
\begin{example}\label{fandr}
Consider a multilattice $\mathcal{M}$ whose partial order is depicted by the diagram in Figure~\ref{fig:rf-sub m-lattice}.
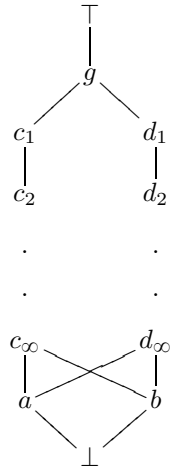
\begin{figure}[!h]
\vspace*{2mm}
	\begin{center}
	$\xymatrix@-1pc@M= 1.5pt{ &\top&  \\
&g\ar@{-}[u] &  \\
c_1\ar@{-}[ur]& &d_1\ar@{-}[ul] \\
c_2\ar@{-}[u]& &d_2\ar@{-}[u] \\
.& &. \\
.& &. \\
c_\infty& &d_\infty \\
a\ar@{-}[u]\ar@{-}[urr]& &b\ar@{-}[u]\ar@{-}[ull] \\
& \bot\ar@{-}[ur]\ar@{-}[ul]& \\
}
$
	\end{center}\vspace*{-3mm}
	\caption{A multilattice, $\mathcal{M}$, with an $r$-sub multillatice, $M\setminus\{d_\infty\}$, that is not a multilattice on its own.} \label{fig:rf-sub m-lattice}
\end{figure}

The subset $A:=M\setminus\{d_{\infty}\}$ is an $r$-sub-multilattice but the restriction of the order on $M$ to $A$ is not a multilattice. Indeed $a\sqcup b=\{c_{\infty}\}$, $a, b\leq d_1$ but there is no $x\in a\sqcup b$ such that $x\leq d_1$.
\end{example}

\begin{proposition}\label{containsM6}
Every pure and bounded multilattice contains a restricted submultilattice isomorphic to $\mathrm{ML_6}$.
\end{proposition}

\begin{proof} Let $(P,\le)$ be a pure and bounded multilattice.
Then, there exists at least two elements $ x, y\in P$ such that $x\sqcap y$ or $ x\sqcup y$ is not a singleton.
\begin{itemize}
    \item If $ x\sqcup y$ is not a singleton, then there are $a, b\in x\sqcup y$ with $a\ne b$. Let $c\in a\sqcup b$ and $z\in x\sqcap y$; then $\{z, x, y, a, b, c\}$ is a restricted submultilattice which is isomorphic to $\mathrm{ML_6}$.
   \item If $x\sqcap y$ is not a singleton, then there are $a, b\in x\sqcap y$ with $a\ne b$. Let $c\in a\sqcap b$ and $z\in x\sqcup y$; then
$\{c, a, b, x, y, z\}$ is an isomorphic copy of $\mathrm{ML_6}$, which is indeed a restricted submultilattice of $P$.
\end{itemize}

\vspace*{-6mm}
\end{proof}

A pure and bounded multilattice does not always contain a full submultilattice isomorphic to $\mathrm{ML_6}$, as the following example shows.

\begin{figure}[!h]
\vspace*{-1mm}
	\begin{center}
$\xymatrix@-1pc@M= 1.5pt{&\top& \\
f\ar@{-}[ur]&& e\ar@{-}[ul] \\
d\ar@{-}[u]\ar@{-}[urr]&& c\ar@{-}[u]\ar@{-}[ull] \\
a\ar@{-}[u]\ar@{-}[urr]&& b\ar@{-}[u]\ar@{-}[ull] \\
& \bot\ar@{-}[ul]\ar@{-}[ur]& }$
\end{center}\vspace*{-3mm}
\caption{A pure bounded multilattice with no full  submultilattice .}\label{fig:pmlatnofsubmlat}\vspace*{-4mm}
\end{figure}
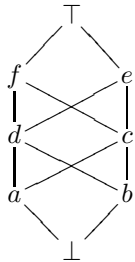

\begin{example}
Consider the multilattice depicted in the Hasse diagram on Figure~\ref{fig:pmlatnofsubmlat}.
All the restricted submultilattices of this multilattice isomorphic to $\mathrm{ML_6}$ are:
\[\{\bot,a, b, c, d, e\},\quad \{\bot,a, b, c, d, f\},\quad \{a, c, d, e, f, \top\} \text{ and } \{b, c, d, e, f, \top\}.
\]
 None of these is a full  submultilattice.
\end{example}

After recalling the concept of complete multilattice, we proceed to introduce residuated multilattices. These will serve as truth degree sets to evaluate objects and attributes.

\section{Some constructions on residuated multilattices}\label{s:TDS}

In order to use residuated multilattices in FCA to evaluate the attributes and objects, we should make sure that there are enough examples. 
To this aim  we will construct new residuated multilattices from others. First we set the notation and terminology for residuated multilattices.

\begin{definition}\cite{JT02}
A structure $\mathcal{A}:= (A,\le, \top, \odot, \rightarrow)$ is called \neu{pocrim} (partially ordered
commutative residuated integral monoid) if  $(A,\leq,\top)$ is a poset with a maximum $\top$ and
$(A, \odot, \top)$ is a commutative monoid such that
\[
(\ddagger)\quad a\odot c\leq b\iff c\leq a\to b, \quad \text{for all } a, b, c\in A.
\]
Any pair $(\odot,\to)$ satisfying $(\ddagger)$ is called a \neu{residuation}, a \neu{residuated couple} or an \neu{adjoint couple} on $(A,\le)$. 	
\end{definition}
%

Let $\mathcal{A}:= (A,\le, \top, \odot \rightarrow)$ be a pocrim and $a, b, c \in A$. The following properties hold:
\begin{itemize}
	\itemsep4pt
	\item[\textbf{P1}] $a\odot b \leq a$ and $a\odot b \leq b$;
	\item[\textbf{P2}]
	$\left\{\begin{array}{rcccl}
	a\odot(a\rightarrow b)&\leq& a &\leq& b\rightarrow (a\odot b) \\
	b\odot(a\rightarrow b)&\leq& b &\leq& a\rightarrow (a\odot b);
	\end{array}\right.$
	\item[\textbf{P3}] $a\leq b$ $\Leftrightarrow$ $a\rightarrow b=\top$;
    \item[\textbf{P4}] $a\leq b\implies\left\{\begin{array}{rcl}
	a\odot c &\le& b \odot c,\\
	c\rightarrow a &\le& c\rightarrow b,\\
	b\rightarrow c &\le & a\rightarrow	c;
	\end{array}  \right. $	
    \item[\textbf{P5}] $a\rightarrow (b\rightarrow c)\ =\ b\rightarrow (a\rightarrow c)$;
    \item[\textbf{P6}] $(a\rightarrow b)\odot (b\rightarrow c)\ \leq\ a\rightarrow c$;
    \item[\textbf{P7}]  $\left\{\begin{array}{rcl}
    a\rightarrow b &\leq & (a\odot c)\rightarrow (b\odot c) \\
	a\rightarrow b &\leq & (c\rightarrow a)\rightarrow (c\rightarrow b) \\
	a\rightarrow b &\leq & (b\rightarrow c)\rightarrow (a\rightarrow c)
	\end{array}\right.$
\end{itemize}

\noindent A pocrim is \neu{bounded} if it also has a lower bound, usually denoted by $\bot$.
A nonempty subset $F$ of a pocrim $A$ is called a \neu{deductive system} (\neu{$ds$}, for short) if
\begin{itemize}
\itemsep=0.85pt
    \item[(ds-1)]  $\top \in F$ and
    \item[(ds-2)] if $x, x\to y\in F$, then $y\in F$ for all $x, y\in A$;
\end{itemize}
    or equivalently
\begin{itemize}
\itemsep=0.85pt
    \item[(i)] $x\odot y\in F$ for all $x, y\in F$ and
    \item[(ii)] if $x\leq y$ and $x\in F$, then $y\in F$ for all $x, y\in A$.
\end{itemize}

Residuated multilattices were introduced in \cite{CCMO14} as follows:

\begin{definition}\label{RML}\cite{CCMO14}
A \neu{residuated multilattice} (write $\mathcal{RML}$ for short)  is a pocrim whose underlying poset is a multilattice.
A \neu{complete residuated multilattice}  is a pocrim whose underlying poset is a complete multilattice.
\end{definition}
A \neu{residuated lattice} is a pocrim whose underlying poset is a lattice.
From Definition \ref{RML}, it is clear that every $\mathcal{RML}$ is a full multilattice.
Let $\mathcal{A}$ be a bounded $\mathcal{RML}$. For simplicity we will write   $x^\ast=x\rightarrow \bot$ and
     $X^\ast=\{x^\ast,~~x\in X\}$ for all  $x\in M$ and  $X\subseteq M$.
A map $h: M\rightarrow P$ between residuated multilattices is said to be a \textbf{homomorphism} if h is a multilattice homomorphism that satisfies $h(a\odot b)=h(a)\odot h(b)$ and $h(a\rightarrow b)=h(a)\rightarrow h(b)$ for all $a$, $b\in M$.

\begin{definition}\label{def:f-sub-m-lattice}
Let $\mathcal{A}$ be an $\mathcal{RML}$ and $X$ a subset of $A$. We say that $X$ is a  \neu{full residuated sub-multilattice} (or \neu{$f$-Sub-$\mathcal{RML}$} for short) if the following conditions hold.
\begin{itemize}
\itemsep=0.85pt
\item[S1.] $\top \in X$
\item[S2.] For every $x, y\in X$, $x\odot y\in X$, $x\to y\in X$.
\item[S3.] $X$ is a $f$-Sub-multilattice.
\end{itemize}
If we replace S3 in Definition~\ref{def:f-sub-m-lattice} by "$X$ is a restricted sub-multilattice", then we obtain the definition of a \neu{restricted residuated sub-multilattice} (or \neu{$r$-Sub-$\mathcal{RML}$} for short).
\end{definition}

For convenience we summarize the main properties of residuated multilattices needed throughout this paper. They can be found or derived from some properties in \cite{CCMO14}.

\begin{proposition}\cite{CCMO14}\label{conditions}
In an $\mathcal{RML}$ $\mathcal{A}$, the following conditions
hold, for all $x, y, z \in A$:
\begin{center}
\begin{enumerate}
\itemsep=0.85pt
  \item[\textbf{M1}] $x\odot y, x\odot(x\rightarrow y)\in ~~~~ \downarrow(x\sqcap y)$;
  \item[\textbf{M2}] $(x\odot y)\sqcup(x\odot z)\subseteq x\odot(y\sqcup z)$;
  \item[\textbf{M3}] $(x\sqcap y)\rightarrow z\subseteq \uparrow[(x\rightarrow z)\sqcup(y\rightarrow z)]$;
  \item[\textbf{M4}] $(x\sqcup y)\rightarrow z\subseteq \downarrow[(x\rightarrow z)\sqcap(y\rightarrow z)]$;
  \item[\textbf{M5}] $(x\rightarrow z)\sqcap(y\rightarrow z)\subseteq (x\sqcup y)\rightarrow z$;
  \item[\textbf{M6}] $x\rightarrow y=\max\{ (x\sqcup y)\rightarrow y \}=\max\{ x\rightarrow (x\sqcap y)\}$;
      \item[\textbf{M7}] $x\leq x^{\ast\ast}$, $x^{\ast}= x^{\ast\ast\ast}$, $x^{\ast\ast}\rightarrow y^{\ast\ast}=y^{\ast}\rightarrow x^{\ast}$;
   \item[\textbf{M8}] $(x\sqcap y)^{\ast}\subseteq\uparrow(x^{\ast} \sqcup y^{\ast}) $;
 \item[\textbf{M9}] $(x\sqcup y)^{\ast}\subseteq \downarrow (x^{\ast}\sqcap y^{\ast})$;
      \item[\textbf{M10}]
$(x^{\ast} \sqcap y^{\ast})\subseteq (x\sqcup y)^{\ast}$.
 \end{enumerate}
 \end{center}
\end{proposition}
Given $\mathcal{A}$, an $\mathcal{RML}$, a nonempty subset $F$ of $A$ is called a filter if $F$ is a $ds$ satisfying: for all $x, y\in A$, if $x\to y\in F$, then $x\sqcup y\to y\subseteq F$ and $x\to x\sqcap y\subseteq F$.

\begin{definition}\cite{CCMO14}
Let  $\mathcal{A}$  be a residuated multilattice. A ds $F$ is said to be \neu{consistent} if for
all $x, y, z \in A$  the following conditions hold:
\begin{itemize}
\itemsep=0.85pt
  \item[CF-1] If $x \rightarrow z, y \rightarrow z \in F$, then $(x \sqcup y) \rightarrow z \subseteq F$.
  \item[CF-2] If $z \rightarrow x, z\rightarrow y \in F$, then $z \rightarrow (x \sqcap y)\subseteq F$.
\end{itemize}
\end{definition}
It is known that consistent ds are filters \cite[Prop. 55]{CCMO14}.
Given a multilattice $\mathcal{A}$ , a relation $\mathcal{R}$ on $A$ and $X, Y\subseteq A$, we write $X\widehat{\mathcal{R}} Y$ if for every $x\in X$, there exists $y\in Y$ such that $x\mathcal{R} y$, and for every $y\in Y$, there exists $x\in X$ such that $x\mathcal{R} y$.
A congruence on $\mathcal{A}$ is an equivalence relation $\mathcal{R}$ on $A$ such that $x\mathcal{R} y$ implies $x\sqcup z\widehat{\mathcal{R}}y\sqcup z$, $x\sqcap z\widehat{\mathcal{R}}y\sqcap z$, $x\odot z\mathcal{R} y\odot z$, $x\to z\mathcal{R} y\to z$, $z\to x\mathcal{R} z\to y$, for all $x, y, z\in A$.\\
If $\mathcal{R}$ is a congruence on an $\mathcal{RML}$ $\mathcal{A}$, then the quotient set $\mathcal{A}/\mathcal{R}=\{[x]:x\in A\}$ is an $\mathcal{RML}$, with the top element $[\top]$ and the operations defined for $x, y\in A$ by:
\begin{align*}
    [x]\,\sqcup [y] &=\{[t]:t\in x\sqcup y\} \\
    [x]\,\sqcap [y] &=\{[t]:t\in x\sqcap y\} \\
    [x]\,\odot [y] &=[x\odot y] \\
    [x]\to\! [y] &=[x\to y]
\end{align*}

Given a filter $F$ of an $\mathcal{RML}$ $\mathcal{A}$, the relation $\mathcal{R}_F$ defined on $\mathcal{A}$ by $x\mathcal{R}_F y$ if and only if $x\to y, y\to x\in F$ is a congruence relation on $\mathcal{A}$ \cite[Thm. 53]{CCMO14}. The quotient set $\mathcal{A}/\mathcal{R}_F$ will be denoted by $\mathcal{A}/F$.
It is known that a filter $F$ is consistent if and only if $\mathcal{A}/F$ is a residuated lattice \cite[Corollary 58]{CCMO14}.

\medskip
The second step consists in identifying the smallest pure $\mathcal{RML}$. Our aim is to prove by contradiction that there is no pocrim structure on $\mathrm{ML_6}$ and to give an example of pocrim structure on a pure and complete multilattice with $7$ elements.

\begin{lemma}\label{fleche}
	Suppose that the multilattice $\mathrm{ML_6}$ is endowed with a pocrim structure
$(\mathrm{ML}_6, \leq, \odot, \to, \bot, \top)$. Then the implication $\to$ should be given by Table~\ref{tab:tableML_6}.

\begin{table}[h]
\caption{A candidate $\to$ on $\mathrm{ML}_6$} \label{tab:tableML_6}\vspace*{-2mm}
\centering
 \scalebox{0.9}{
 	\begin{tabular}{|c|c|c|c|c|c|c|}
		\hline
		$\rightarrow$  & $\bot$ & a & b & c & d & $\top$ \\
		\hline
		$\bot$& $\top$ & $\top$ & $\top$ & $\top$ & $\top$ & $\top$ \\
		\hline
		a & b & $\top$ & b & $\top$ &$\top$ & $\top$ \\
		\hline
		b & a & a & $\top$ &  $\top$ &$\top$ & $\top$ \\
		\hline
		c & $\bot$ & a & b & $\top$ & d & $\top$\\
		\hline
		d & $\bot$ & a & b & c  & $\top$ & $\top$\\
		\hline
		$\top$ &$\bot$ & a & b & c & d & $\top$ \\
		\hline
	\end{tabular} }
\end{table}\vspace*{-2mm}
\end{lemma}

\begin{proof}
	We shall use repeatedly the fact that in any pocrim $y\leq x\to y$ and $x\to (y\to z)=y\to (x\to z)$ hold.
	\begin{itemize}
		\item[(i)] Since $c\leq d\to c\ne \top$, then $d\to c=c$. Similarly $c\to d=d$.
		\item[(ii)] Note that since $b\leq c\to b$ and $b\leq d\to b$, then $c\to b\in \{b, c, d\}$ and
$d\to b\in \{b, c, d\}$. We narrow down the possibilities further. Note that
$c\to b\leq c\to d=d$, hence $c\to b\ne c$. Likewise, $d\to b\ne d$. Thus $c\to b\in \{b, d\}$ and
$d\to b\in \{b, c\}$. We show that $c\to b=b$ and $d\to b=b$ by showing that the other combinations lead to
contradiction.
		\begin{itemize}
			\item If $c\to b=b$ and $d\to b=c$, then

 $c=d\to b=d\to (c\to b)=c\to (d\to b)=c\to c=\top$, which is a contradiction.
			\item If $c\to b=d$ and $d\to b=b$, then

$\top=d\to d=d\to (c\to b)=c\to (d\to b)=c\to b$, which implies that $c\leq b$. This is impossible.
			\item If $c\to b=d$ and $d\to b=c$, then

$\top =a\to c=a\to (d\to b)=d\to (a\to b)$. So $d\leq a\to b\ne \top$, hence $d=a\to b$. Now,
$\top=a\to d=a\to (c\to b)=c\to (a\to b)=c\to d$, which is again a contradiction.
\end{itemize}
Since $a$ and $b$ play symmetrical roles, we deduce $c\to a=a$ and $d\to a=a$.
		\item[(iii)] Since $b\leq a\to b$, then $a\to b \in \{b, c, d\}$. Again, we show that
$a\to b \in \{c, d\}$ is impossible. Indeed, if $a\to b=c$, then

$\top=c\to c=c\to (a\to b)=a\to (c\to b)=a\to b$, which is impossible. Similarly, if $a\to b=d$, then
$\top=d\to d=d\to (a\to b)=a\to (d\to b)=a\to b$, which is the same contradiction. Hence, $a\to b=b$. The proof
that $b\to a=a$ is analogous.
		\item[(iv)] It remains to show that $a\to \bot=b, b\to \bot =a$ and $c\to \bot=d\to \bot=\bot$.\\
		Note that $a\odot b\leq a, b$, so $a\odot b=\bot$. Thus, $a\leq b\to \bot$ and $b\leq a\to \bot$. On the
other hand, $b\to \bot\leq b\to a=a$. Hence, $b\to \bot =a$ and a similar argument shows that $a\to \bot =b$.
Finally, observe that $a, b\leq c$, so $c\to \bot\leq a\to \bot, b\to \bot$. Hence, $c\to \bot\leq b, a$ and
consequently $c\to \bot=\bot$. The verification that $d\to \bot =\bot$ is similar.
	\end{itemize}
	We have justified all the entries of Table~\ref{tab:tableML_6}.
\end{proof}

\begin{proposition}\label{noprml}
	There does not exist a residuated multilattice structure on $\mathrm{ML_6}$ extending its existing partial order.
\end{proposition}

\begin{proof}
	By contradiction suppose that there exist $\odot$ and $\to$ such that
$(\mathrm{ML}_6,\leq, \odot, \to, \bot, \top)$ is a residuated multilattice. Then by Lemma \ref{fleche}, $\to$
is given by table 1 above.\medskip

Since $a\odot a\leq a$, we have $a\odot a\in \{\bot, a\}$.

If $a\odot a=\bot$, then $a\odot a\leq b$ and $a\leq a\to b=b$ (Table 1), which is a contradiction. \medskip

Suppose $a\odot a=a$. Since $a\leq c$, we have $a=a\odot a\leq a\odot c\leq a$. Hence, $a\odot c=a$. From
$a\leq d$, we have $a=c\odot a\leq c\odot d\leq c, d$. Because $c$ and $d$ are incomparable, $c\odot d=a$. It
follows that $d\leq c\to a=a$ (Table 1), which is again a contradiction.
\end{proof}

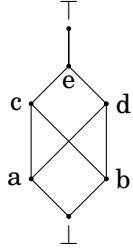
\begin{figure}[!ht]
	\begin{center}
		\begin{tikzpicture}[scale=0.5]
		\draw[fill] (2,0) circle (0.05);  
		\draw[fill] (1,1) circle (0.05);  
		\draw[fill] (3,1) circle (0.05);  
		\draw[fill] (3,3) circle (0.05);  
		\draw[fill] (1,3) circle (0.05);  
		\draw[fill] (2,4) circle (0.05);  
		\draw[fill] (2,5) circle (0.05);  
		\draw (2,0) -- (1,1) -- (1,3) -- (2,4) -- (2,5) -- (2,4) --(3,3)--(3,1)--(2,0);
		\draw (1,1) -- (3,3);
		\draw (1,3) -- (3,1);
		\draw (2,0) node[below]{$\bot$};  
		\draw (1,1) node[left]{a};  
		\draw (3,1) node[right]{b};  
		\draw (3,3) node[right]{d};  
		\draw (1,3) node[left]{c};  
		\draw (2,4) node[below]{e};  
		\draw (2,5) node[above]{$\top$};  
		\end{tikzpicture}
	\end{center}\vspace*{-5mm}
	\caption{The Hasse diagram of a complete residuated multilattice with seven element: $\mathrm{RML}_7$.}\label{fig:rcmlat}
\end{figure}

\noindent Combining Proposition~\ref{containsM6} and Proposition~\ref{noprml} we obtain that every bounded pure
$\mathcal{RML}$ has  at least seven elements.
The next example shows that there is indeed a bounded
pure $\mathcal{RML}$ with seven elements. We will denote it by $\mathrm{RML}_7$. It Hasse diagram is given by Figure~\ref{fig:rcmlat} and its operations $\odot$ and $\rightarrow$ defined in the following tables.

\begin{table}[h]\small
\caption{Tables of $\odot$ and $\rightarrow$ on $\mathrm{RML}_7$}\label{tableML_7}\vspace*{-1mm}
\centering
\begin{subtable}{0.4\linewidth}
\hspace*{-7mm}\begin{tabular}{|c|c|c|c|c|c|c|c|}
		\hline
		$\odot$ & $\bot$ &    a   &    b   &    c   &    d   &    e   & $\top$ \\
		\hline
		$\bot$  & $\bot$ & $\bot$ & $\bot$ & $\bot$ & $\bot$ & $\bot$ & $\bot$ \\
		\hline
		   a    & $\bot$ &    a   & $\bot$ &    a   &    a   &    a   &  a \\
		\hline
		   b    & $\bot$ & $\bot$ & $\bot$ & $\bot$ & $\bot$ & $\bot$ &  b \\
		\hline
		   c    & $\bot$ &    a   & $\bot$ &    a   &    a   &     a  &   c \\
		\hline
		   d    & $\bot$ &    a   & $\bot$ &    a   &    a   &     a  &  d \\
		\hline
		   e    & $\bot$ &    a   & $\bot$ &    a   &    a   &    a   &  e \\
		\hline
		$\top$  & $\bot$ &    a   &    b   &    c   &    d   &    e   & $\top$ \\
		\hline
	\end{tabular}
\caption{Table of $\odot$ on $\mathrm{RML}_7$}\label{table 1}
\end{subtable}
             \hspace*{10mm}
\begin{subtable}{0.4\linewidth}
	\begin{tabular}{|c|c|c|c|c|c|c|c|}
		\hline
		$\rightarrow$  & $\bot$ &    a   &    b   &    c   &    d   & e & $\top$ \\
		\hline
		     $\bot$    & $\top$ & $\top$ & $\top$ & $\top$ & $\top$ & $\top$ & $\top$ \\
		\hline
		        a      &    b   & $\top$ &    b   & $\top$ & $\top$ & $\top$ & $\top$ \\
		\hline
		        b      &    e   &    e   & $\top$ & $\top$ & $\top$ & $\top$ & $\top$ \\
		\hline
		        c      &    b   &    e   &    b   & $\top$ &    e   & $\top$ & $\top$\\
		\hline
		        d      &    b   &    e   &    b   &    e   & $\top$ & $\top$ & $\top$\\
		\hline
		        e      &    b   &    e   &    b   &    e   &    e   & $\top$ & $\top$\\
		\hline
		     $\top$    & $\bot$ &    a   &    b   &    c   &    d   &    e   & $\top$ \\
		\hline
	\end{tabular}
\caption{Table of $\to$ on $\mathrm{RML}_7$.}\label{table 2}
\end{subtable}\vspace*{-2mm}
\end{table}

\medskip
Now that we have set up a smallest (in term of the number of elements) pure $\mathcal{RML}$, we proceed to apply the ordinal sum construction to $\mathcal{RML}$s and obtain other (pure) $\mathcal{RML}$s.
First, we recall the construction of the ordinal sum of pocrims (see \cite{AU21, JM09, Bm04}). We take two pocrims $(A,\leq_A, \odot_A, \to_A, \top_A)$
and $(B, \leq_B, \odot_B,\to_B, \top_B)$ and define the set
$C:=(A\amalg B) / \{\top_A\cong\top_B\}$ (that is the disjoint union of $A$ and $B$ with $\top_A$ and $\top_B$
identified), and a relation $\leq $ on $C$ by $x\leq y$ if $x, y\in A$ and $x\leq_A y$, or $x, y\in B$ and $x\leq_B y$
or $x\in A\setminus \{\top_A\}$ and $y\in B$.

\medskip
The operation $\odot$ and $\to$ on $C$ are defined by:
\begin{align*}
	x\odot y=
	\begin{cases}
	x\odot_A y &\text{ if }\ x, y\in A \\
	x\odot_B y &\text{ if }\ x, y\in B \\
	x &\text{ if }\ x\in A\setminus \{\top_A\}\ \text{ and }\ y\in B  \\
	y &\text{ if }\  x\in B \text{ and }\ y\in A\setminus \{\top_A\}
	\end{cases}	\vspace*{-2mm}
\end{align*}
\begin{align*}
	x\rightarrow y=
	\begin{cases}
		x\to_A y &\text{ if }\ x, y\in A  \\
	x\to_B y &\text{ if }\ x, y\in B  \\
	\top_B & \text{ if }\ x\in A\setminus \{\top_A\}\ \text{ and }\ y\in B  \\
	y & \text{ if }\ y\in A\setminus \{\top_A\}\ \text{ and }\  x\in B
	\end{cases}
\end{align*}

Then $(C, \leq, \odot, \to, \bot_A, \top_B)$ is a pocrim, called the ordinal sum of $A$ and $B$ and denoted
by $A\oplus B$. Note that in $A\oplus B$, $A\setminus\{\top_A\}$ is below every element of $B$. Our goal is to apply this construction to create new $\mathcal{RML}$ from old ones.

\begin{proposition}\label{ordsum}
Let $\mathcal{M}$ and $\mathcal{N}$ be bounded $\mathcal{RML}$s. Then
\begin{itemize}
\item[(1)] The ordinal sum $M\oplus N$ (as pocrims) is an $\mathcal{RML}$.
\item[(2)] $M$ is an f-Sub-$\mathcal{RML}$ of $M\oplus N$ if and only if $M$ has a unique coatom.
\item[(3)] $N$ is a filter of $M\oplus N$, and in particular an $f$-Sub-$\mathcal{RML}$.
\item[(4)] The quotient $\mathcal{RML}$ $M\oplus N/N$ is canonically isomorphic to $M$.
\item[(5)] $N$ is a consistent filter of $M\oplus N$ if and only if $\mathcal{M}$ is a residuated lattice.
\end{itemize}
\end{proposition}
\begin{proof}
Let  $\mathcal{M}$, $\mathcal{N}$ be two $\mathcal{RML}$s.

\begin{itemize}
\item[(1)]
    We know that the ordinal sum of $\mathcal{M}$ and $\mathcal{N}$ (as pocrims) is again a pocrim. In addition, it is clear that with respect to the ordinal sum's order, $\sqcup_{M\oplus N}$ and $\sqcap_{M\oplus N}$ are given by:
\begin{align*}
	x\sqcup_{M\oplus N}y &=
	\begin{cases}
	x\sqcup_{M}y &\text{ if }\ x, y\in M \text{ and } x\sqcup_{M} y\ne \top_{M} \\
	x\sqcup_{N} y &\text{ if }\ x, y\in N \\
	y &\text{ if }\ x\in M\setminus \{\top_M\}\ \text{ and }\ y\in N  \\
	\bot_{\mathcal{N}} &\text{ if }\  x,y\in M\setminus\{\top_M\} \text{ and }\ x\sqcup y =\top_M
	\end{cases}	 \\[5pt]
	x\sqcap_{M\oplus N}y &=
	\begin{cases}
	x\sqcap_{M}y &\text{ if }\ x, y\in M \\
	x\sqcap_{N} y &\text{ if }\ x, y\in N \\
	x &\text{ if }\ x\in M \text{ and }\ y\in N.
	\end{cases}
\end{align*}

It remains to show that for every $x, y, a, b\in M\oplus N$, with $x, y\leq a$ and $b\leq x, y$, there exists $u\in x\sqcup_{M\oplus N}y$ and $v\in x\sqcap_{M\oplus N}y$ such that $u\leq a$ and $b\leq v$.
This is easily verified by considering the cases $a\in M\setminus \{\top_{M}\}$, $a\in N$ and $x, y\in M\setminus \{\top_{M}\}$ or $x, y\in N$ or $x\in M\setminus \{\top_{M}\}$ and $y\in N$ as necessary.
Therefore $M\oplus N$ is an $\mathcal{RML}$ as claimed.
\item[(2)] Assume there exists $x_0\ne y_0$ coatoms in $M$. Then $x_0\sqcup_{M\oplus N} y_0=\bot_{\mathcal{N}}\notin M$. Thus $M$ is not an f-Sub-$\mathcal{RML}$ of $M\oplus N$.

Conversely suppose that $M$ has a unique coatom $a$. Then for all $x, y\in M\setminus \{\top_{M}\}$, $x\sqcup_{M\oplus N}y=x\sqcup_{M} y\subseteq \downarrow a\subseteq M$. If $x=\top_{M}$ or $y=\top_{M}$, then $x\sqcup_{M\oplus N}y=\top_{M}\in M$. It is clear from the description of $\sqcap_{M\oplus N}$ above that $x\sqcap_{M\oplus N}y\subseteq M$ for all $x, y\in M$. In addition it follows from the definition of the ordinal sum of pocrims that $x\to y, x\odot y\in M$ for all $x, y\in M$. Therefore, $M$ is an f-Sub-$\mathcal{RML}$ of $M\oplus N$.

\item[(3)]To show that $N$ is a filter of $M\oplus N$, we first show that $N$ is a deductive system. It follows from the definition of $\odot$ in the ordinal sum that $N$ is $\odot$-closed, and from the definition of the order that whenever $x\leq y$ with $x\in N$, then $y\in N$. So, $N$ is a deductive system.
Moreover, let $x, y\in M\oplus N$ such that $x\rightarrow y\in N$. By definition of $\to$ in the ordinal sum, either $x, y\in N$ and $x\to y=x\rightarrow_{N} y$  or $x\to y=\top_{N}$. If $x\rightarrow y=x\rightarrow_{N} y$ with $x, y\in N$ then from the description of $\sqcup_{M\oplus N}$ above and that of $\to$ in the ordinal sum, it follows that $x\sqcup_{M\oplus N} y\to y\subseteq N$ and $x\to x \sqcap_{M\oplus N}y\subseteq N$. If $x\to y=\top_{N}$, then $x\leq y$. Thus, $x\sqcup y\to y=y\to y=\top_{N}\in N$  and $x\to x\sqcap y=x\to x=\top_{N}$.\\
Thus $N$ is a filter of $M\oplus N$.
\item[(4)] From the definition of $\to$ in the ordinal sum, it is clear that for every $x\in M\oplus N$, $[x]_{N}=\{x\}$ if  $x\in M\setminus\{\top_{M}\}$ and $[x]_{N}=N$ otherwise.
Now consider $f:M\oplus N/N\to M$ defined by:
   \begin{eqnarray*}
          && f([x]_{N})= \left\{
                     \begin{array}{ll}
                        x \; \; \; \;\; \; \; \;~~~~\text{if} ~~ x\in M\setminus\{\top_{M}\} & \hbox{} \\
                        \top_{M} \; \; \; \;\; \; \; \;~~~~\text{otherwise}  & \hbox{}
                     \end{array}
                   \right.
    \end{eqnarray*}
An elementary but lengthy argument shows that $f$ is a well-defined isomorphism of $\mathcal{RML}$s.
\item[(5)] We recall that a filter $F$ of an $\mathcal{RML}$ $\mathcal{A}$ is consistent if and only if $M/F$ is a residuated lattice (see for e.g., \cite[Corollary 57]{CCMO14}). Therefore, the result is a consequence of (4).
\end{itemize}

\vspace*{-9mm}
\end{proof}
As an application of this construction, we address the comparison of maximal and consistent filters.

\begin{remark}
We would like to point out that there is no general comparison between consistent filters and maximal filters. Indeed, since every filter of a residuated lattice is consistent, then it follows that a consistent filter needs not be maximal. Conversely, if one considers the ordinal sum of two copies of $RML_7$, then $RML_7$ is a maximal filter of $RML_7\oplus RML_7$ that is not consistent since $RML_7\oplus RML_7/RML_7\cong RML_7$, which is not a residuated lattice.
\end{remark}

\begin{figure}[!b]
    \centering
$\xymatrix@-1pc@M= 1.5pt{ & \top & \\
& c\ar@{-}[u]& \\
a\ar@{-}[ur]& &b\ar@{-}[ul] \\
& \bot\ar@{-}[ur]\ar@{-}[ul]& \\
}
$
\caption{Hasse diagram of $RL_5$}   \label{fig:RL5}
\end{figure}
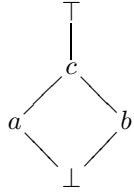

Next, we apply the ordinal sum to construct a new example of an $\mathcal{RML}$.
\begin{example}\label{ex sum ord M5} We wish to construct a concrete $\mathcal{RML}$ as the ordinal sum of  $RML_7$  and the residuated lattice $RL_5$ depicted in the Figure~\ref{fig:RL5}:
where the multiplication and implication are defined by $x\odot y=x\wedge y$, $x\to y=\top$ if  $x\leq y$,  $x\to y=y$ if $x=\top$ or $x=c$ and $y<c$, $a\to \bot=a\to b=b, b\to \bot= b\to a=a$.
\eject

\noindent For simplicity, we rename elements in the ordinal sum as:
$$a_0, a_1, a_2, a_3, b_0, b_1, b_2, b_3, b_4, b_5, b_6$$
and set $A:=\{a_0, a_1, a_2, a_3\}$,  $B:=\{b_0, b_1, b_2, b_3, b_4, b_5, b_6\}$ and  $C:=\{b_2, b_3, b_4, b_5 \}$. \\We shall obtain an $\mathcal{RML}$ whose  Hasse diagram is depicted as:
$$
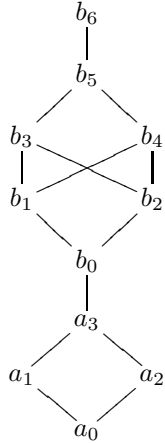
\xymatrix@-1pc@M= 1.5pt{ &b_6&  \\
&b_5\ar@{-}[u]&  \\
b_3\ar@{-}[ur]& &b_4\ar@{-}[ul] \\
b_1\ar@{-}[u]\ar@{-}[urr]& &b_2\ar@{-}[u]\ar@{-}[ull] \\
& b_0\ar@{-}[ur]\ar@{-}[ul]& \\
& a_3\ar@{-}[u]& \\
a_1\ar@{-}[ur]& & a_2\ar@{-}[ul] \\
& a_0\ar@{-}[ur]\ar@{-}[ul]& \\
} 
$$
Moreover, the multiplication and implication are given by:
\begin{equation*}
    x\odot y=\left\{
       \begin{array}{ll}
         x\wedge y\; \; \; \; \; \;~~~~~\text{if} ~~ x, y\in A  & \hbox{} \\
b_0 \; \; \; \; \; \; \; \; \; \; \;~~~~~\text{if} ~~ x\in \{b_0, b_1\} ~~\text{and}~~ y\in B\setminus\{b_6\} ~~\text{or}~~~~ y\in \{b_0, b_1\}~~\text{and}~~ x\in B\setminus\{b_6\}& \hbox{} \\
x \; \; \; \; \; \; \; \;\; \; \; \;~~~~~\text{if} ~~ y=b_6~~\text{or}~~x\in A ~and~ y\in B\setminus\{b_6\}  & \hbox{} \\
y \; \; \; \; \; \; \;\; \; \; \; \;~~~~~\text{if} ~~ x=b_6~~\text{or}~~y\in A ~and~ x\in B\setminus\{b_6\}  & \hbox{} \\
b_2~~~~~~~~~~~~~~~~\text{otherwise}
       \end{array}
     \right.
\end{equation*}
\begin{equation*}
    x\rightarrow y=\left\{
       \begin{array}{ll}
         b_6\; \; \; \; \; \; \; \; \; \; \; \;~~~\text{if} ~~x\leq y\\
 y \; \; \; \; \; \; \; \; \; \; \; \; \;~~~\text{if}~~ x=b_6~~\text{or}~~~ x =a_3 ~~\text{and}~~ y\in\{a_0, a_1, a_2\}~~~ \text{or}~~~ x\in B\setminus\{b_6\}~~\text{and}~~ y\in A\\
 a_1\; \; \; \; \; \; \; \; \; \; \; \;~~~\text{if}~~ x=a_2 ~~\text{and} ~~ y\in\{a_0, a_1\}\\
 a_2\; \; \; \; \; \; \; \; \; \; \; \;~~~\text{if}~~ x=a_1 ~~\text{and} ~~ y\in\{a_0, a_2\}& \hbox{} \\
b_1 \; \; \; \; \; \; \; \; \; \; \; \;~~~\text{if} ~~ x\in  C  ~and~ x\in\{b_0, b_1\}& \hbox{}\\
b_5\; \; \; \; \; \; \; \; \; \; \; \;~~~\text{otherwise}
       \end{array}
     \right.
\end{equation*}
\end{example}

We are now going to look at the structure on the set of all mappings from a non-empty set to a residuated multilattice.
Given a residuated multilattice $\mathcal{A}$ and a nonempty set $X$, we will denote by $A^X$ the set of all
mappings from $X$ to $A$. An order
on $A^X$  as well as the operations $\odot$, $\to$ are defined pointwise: i.e., for $f_1,f_2\in A^X$ we have
\begin{itemize}
\itemsep=0.9pt
	\item $f_1\leq f_2:\iff f_1(x)\leq f_2(x)$ for all $x\in X$,
	\item $\left(f_1\odot f_2\right)(x):=f_1(x)\odot f_2(x)$, for any $x\in X$,
	\item $\left(f_1\to f_2\right)(x):= f_1(x)\to f_2(x)$, for any $x\in X$.
\end{itemize}
\begin{remark}
	Let $\mathcal{A}$ be a complete residuated multilattice and $X$ a nonempty set. It is easy to see that
$\left(A^X,\le,\odot,\to,\top,\bot \right)$ is a complete residuated multilattice.
\end{remark}

Let $(A_1, \leq_1)$ and $(A_2, \leq_2)$ be two complete
multilattices, $G$ and $M$ two non-empty sets (of objects and attributes), and $(\varphi,
\psi)$ a Galois connection between $(A_1^G,\le_1)$ and $(A_2^M,\le_2)$. Recall that a concept is a pair
$(h,f)\in A_1^G\times A_2^M$ such that $\psi(f)=h$ and $\varphi(h)=f$. Concepts are ordered by
\[(h_1, f_1)\le (h_2, f_2):\iff h_1\leq_1 h_2. 
\]

The following theorem has already been shown in \cite{RM12}, we re-enclose it here in order to achieve greater consistency of the material presented and used.
\begin{theorem}\label{thm:1}\cite{RM12}
Let $(A_i, \leq_i)$, $i=1,2$ be two complete
multilattices, $G$ and $M$ be two non-empty sets  and $(\varphi,\psi)$ be a Galois connection between
$(A_1^G,\le_1)$ and $(A_2^M,\le_2)$. Let $H\subseteq A_1^G$, $F\subseteq A_2^M$ and $\mathcal{C}$ the set of
concepts of $A_1^G\times A_2^M$. Then
\begin{itemize}
	\item[\textup{(i)}] $\sqcap\psi(F)\subseteq \psi(\sqcup F)$  and
	$\sqcap\varphi(H)\subseteq \varphi(\sqcup H)$.
	\item[\textup{(ii)}] The poset $(\mathcal{C},\le)$ of all concepts of $A_1^G\times A_2^M$ is a complete
multilattice, with:

$\bigsqcap\limits_{j\in J}(h_j, f_j):=\{(h,\varphi(h)); h\in \sqcap\{h_j\mid j \in J\}\}$ and

$\bigsqcup\limits_{j\in J}(h_j, f_j):=\{(\psi(f), f); f\in \sqcap\{f_j;\mid j \in J\}\}$.
\end{itemize}
\end{theorem}

We are going to show that starting with two complete residuated multilattices $\mathcal{A}_1$ and $\mathcal{A}_2$,
the set of concepts again forms a complete residuated multilattice.

\section{Residuated concept multilattices} \label{s:RCML}

The structures of truth values in fuzzy logic and fuzzy set theory are usually lattices or residuated
lattices (see \cite{GJ67, HP98, HU96}). Particularly in fuzzy formal concept analysis, residuated lattices are used to evaluate the attributes and objects. But the necessity of the use of a more general structure arise in some
examples. In  \cite{RM12} multilattices are used as the underlying set of truth-values in FCA.

\medskip
From now on, $\mathcal{A}_i:=(A_i,\le_i,\top_i,\odot_i,\to_i,\bot_i)$, $i=1,2$ are complete residuated
multilattices, $G$ is a set of objects (to be evaluated in $\mathcal{A}_1$), $M$ is a set of attributes
(to be evaluated in $\mathcal{A}_2$), and $(\varphi,\psi)$ is a Galois connection between $A_1^G$ and $A_2^M$.
Further we  denote by $\mathcal{C}$ the set of concepts,  $\Ext(\mathcal{C})$ the set of extents and
$\Int(\mathcal{C})$ the set of intents. i.e.,
\begin{align*}
	\mathcal{C}&:=\{(h,f)\in A_1^G\times A_2^M ; \varphi(h)=f \text{ and } \psi(f)=h \}, \\
	\Ext(\mathcal{C})&:=\{h\in A_1^G ; (h,\varphi(h))\in \mathcal{C}\} =\{h\in A_1^G ; \psi\varphi(h)=h\},\\
	\Int(\mathcal{C})&:=\{f\in A_2^M ; (\psi(f),f)\in \mathcal{C}\} =\{f\in A_2^M ; \varphi\psi(f)=f\}.
\end{align*}
The operations $\odot_i$ and $\to_i$ are defined componentwise on $A_1^G$ and $A_2^M$. Let $\ddh_i$ and $\uuh_i$
be the constant maps with values $\bot_i$ and $\top_i$: i.e. for $g\in G$ and $m\in M$,
 \[
 \ddh_1(g)=\bot_1,\quad \ddh_2(m)=\bot_2,\quad \uuh_1(g)=\top_1 \quad \text{ and }\quad  \uuh_2(m)=\top_2.
 \]
 For any $h\in A_1^G$ and $f\in A_2^M$ we have
\[
h\odot_1\ddh_1\ =\ \ddh_1,\quad  h\odot_1\uuh_1\ =\ h,\quad
h\to_1\uuh_1\ =\ \uuh_1\quad \text{ and } \]
\[f\odot_2\ddh_2\ =\ \ddh_2,\quad  f\odot_2\uuh_2\ =\ f,\quad f\to_2\uuh_2\ =\ \uuh_2.
\phantom{\quad \text{ and }}\]
We are interested in constructing a residuated couple on $(\mathcal{C},\le)$; We are looking for a suitable
product ($\odot$) and implication ($\to$) such that \[(h_1,f_1)\odot (h_2,f_2)\le (h_3,f_3) \iff (h_1,f_1)\le
(h_2,f_2)\to (h_3,f_3).\]

\begin{lemma} Let $h_1,h_2\in A_1^G$,  $f_1,f_2\in A_2^M$, and $(\varphi,\psi)$ a Galois connection. Then
\begin{enumerate}
	\itemsep=0.9pt
	\item $\psi\varphi(h_1\odot_1 h_2) = \min\{\psi(f) ; f\in A_2^M$   and $h_1\odot_1 h_2\leq_1 \psi(f)\}$
	\item $\varphi\psi(f_1\odot_2 f_2) = \min\{\varphi(h) ; h\in A_1^G \text{ and } f_1\odot_2 f_2\leq_2
\varphi(h)\}$
	\item $\psi\varphi(h_1\to_1 h_2) =\max\{\psi\varphi(h) ; h\in A_1^G \text{ and } h\odot_1  h_1\le_1 h_2 \}$
	\item $\varphi\psi(f_1\to_2 f_2)=\max\{\varphi\psi(f) ; f\in A_2^M \text{ and } f\odot_2 f_1\leq_2 f_2 \}$.
\end{enumerate}
\end{lemma}

\begin{proof}
	 Let $(\varphi,\psi)$ be a Galois connection.
\begin{enumerate}
\itemsep=0.9pt
\item follows from the fact that $\psi\varphi$ is a closure operator on $(A_1^G,\le_1)$ and the elements $\psi(f)$ are $\psi\varphi$-closed.
\item follows from the fact that $\varphi\psi$ is a closure operator on $(A_2^M,\le_2)$ and the elements $\varphi(f)$ are $\varphi\psi$-closed.
\item  Let $h_1,h_2\in A_1^G$. We set $H= \{\psi\varphi(h); h\in A_1^G \text{ and } h\odot_1 h_1\le_1 h_2\}$.
Then
\begin{align*}
H &=\{\psi\varphi(h); h\in A_1^G \text{ and }  h\odot_1 h_1\leq_1  h_2 \}
  =\{\psi\varphi(h); h\in A_1^G \text{ and }  h\le_1 h_1\to_1  h_2 \} \\
&\subseteq\{\psi\varphi(h); h\in A_1^G \text{ and }
\psi\varphi(h)\le_1\psi\varphi(h_1\to_1  h_2) \}, \quad \text{since } \psi\varphi \text{ is isotone}.
\end{align*}
\noindent
Thus $\psi\varphi(h_1\to_1  h_2)$ is an upper bound of $H$. As $h_1\to_1  h_2\le_1 h_1\to_1  h_2$, we get $\psi\varphi(h_1\to_1  h_2)\in H$. Therefore, \[\psi\varphi(h_1\to_1  h_2)=\max\{\psi\varphi(h); h\in A_1^G \text{ and } h\le_1 h_1\to_1  h_2\}.\]
\item The proof is similar to 3.
\end{enumerate}

\vspace*{-8mm}
\end{proof}

\noindent
Let $h_1,h_2\in A_1^G$ and  $f_1,f_2\in A_2^M$. We define the operations $\otimes_1$, $\otimes_2$,
$\rightrightarrows_1$ and $\rightrightarrows_2$ as follows:
\begin{align*}
h_1\otimes_1 h_2:= \psi\varphi(h_1\odot_1 h_2) && h_1\rightrightarrows_1 h_2:= \psi\varphi(h_1\to_1 h_2) \\
f_1\otimes_2 f_2:= \varphi\psi(f_1\odot_2 f_2) &&
f_1\rightrightarrows_2 f_2:= \varphi\psi(f_1\to_2 f_2).
\end{align*}
\noindent
Observe that the operations $\otimes_1$ and $\rightrightarrows_1$ are defined on all pairs $h_1, h_2\in A_1^G$
but their results are in  $\Ext(\mathcal{C})$.
Similarly, the operations $\otimes_2$ and $\rightrightarrows_2$ are defined on all pairs $f_1, f_2\in A_2^M$
but their results are in  $\Int(\mathcal{C})$. We can then infer that  $(h_1\otimes_1 h_2,
\varphi(h_1\otimes_1 h_2))=(\psi\varphi(h_1\odot_1 h_2), \varphi(h_1\odot_1 h_2))$ is a concept and  also that
$(h_1\rightrightarrows_1 h_2, \varphi(h_1\rightrightarrows_1 h_2)) = (\psi\varphi( h_1\to_1 h_2),
\varphi(h_1\to_1 h_2))$ is a concept. Similarly,
$(\psi(f_1\otimes_2 f_2), f_1\otimes_2 f_2)$ and 
$(\psi(f_1\rightrightarrows_2 f_2), f_1\rightrightarrows_2 f_2)$
are concepts.

\medskip
It is obvious that the operations $\otimes_i$, $i\in \{1, 2\}$ are commutative. Now, let us look about their
associativity and the residuated couple needed.

\begin{lemma}\label{lemme 1} Let $(\varphi,\psi)$ be a Galois connection. If $\Ext(\mathcal{C})$ and
$\Int(\mathcal{C})$ are closed under $\to_i$ then $(\otimes_i, \rightrightarrows_i)$, $i=1,2$
are residuated couples and $\otimes_1$ and $\otimes_2$ are associative.
\end{lemma}

\begin{proof} Let $h_1, h_2, h_3 \in A_1^G$. 
\begin{align*}
h_1\otimes_1 h_2\leq_1 h_3 &\iff \psi\varphi (h_1\odot h_2) \leq_1 h_3\\
&\implies h_1\odot h_2 \le_1 h_3 \\
&\iff h_1\le_1 h_2 \to_1 h_3 \le_1 \psi\varphi(h_2 \to_1 h_3) = h_2\rightrightarrows_2 h_3.
\end{align*}
Henceforth, 
$h_1\otimes_1 h_2\leq_1 h_3\implies h_1\leq_1 h_2\rightrightarrows_1 h_3$.

\medskip
The converse holds if $h_3$ and $h_2\to_1 h_3$ are closed. In fact,
\begin{align*}
h_1\leq_1 h_2\rightrightarrows_1 h_3 &\iff h_1\leq_1 \psi\varphi(h_2\to_1 h_3)\\ 
&\iff h_1\leq_1 h_2\to_1 h_3,\quad \text{assuming $h_2\to_1 h_3$ closed}\\
&\iff h_1\odot_1 h_2\leq_1 h_3 \\
&\iff \psi\varphi(h_1\odot_1 h_2)\leq_1 h_3, \quad \text{ assuming $h_3$ closed} \\
&\iff h_1\otimes_1 h_2\leq_1 h_3.
\end{align*}
Thus $(\otimes_1, \rightrightarrows_1)$ is an adjoint couple on $\Ext(\mathcal{C})$, if $\Ext(\mathcal{C})$
is closed under $\to_1$.

\medskip
Similarly, we can proved that $(\otimes_2, \rightrightarrows_2)$ is an adjoint couple on $\Int(\mathcal{C})$,
if it is closed under $\to_2$.
We still need to check the associativity of $\otimes_1$ and $\otimes_2$. Let $h_1,h_2,h_3$ in $A_1^G$. Then
\begin{align*}
(h_1\otimes_1 h_2)\otimes_1 h_3 &= \psi\varphi((h_1\otimes_1 h_2)\odot_1 h_3) \\ &=\min\{\psi(f);
 f\in A_2^M \text{ and } (h_1\otimes_1
h_2)\odot_1 h_3\leq_1 \psi(f)\},
\end{align*}\vspace*{-1mm}
and\vspace*{-1mm}
\begin{align*}
h_1\otimes_1 (h_2\otimes_1 h_3) &= \psi\varphi(h_1\odot_1 (h_2\otimes_1 h_3)) \\ &=\min\{\psi(f);
 f\in A_2^M \text{ and } h_1\odot_1 (h_2\otimes_1 h_3)\leq_1 \psi(f)\}.
\end{align*}
It is enough to prove that the sets
\[H_1=\{\psi(f); f\in A_2^M \text{ and } (h_1\otimes_1
h_2)\odot_1 h_3\leq_1 \psi(f)\}\]
and
\[
H_2=\{\psi(f); f\in A_2^M \text{ and } h_1\odot_1 (h_2\otimes_1 h_3)\leq_1 \psi(f)\}
\]
are equal. Let $f\in A_2^M$ such that $(h_1\otimes_1 h_2)\odot_1 h_3\leq_1\psi(f)$.
\begin{align*}
(h_1\otimes_1 h_2)\odot_1 h_3\leq_1 \psi(f) &\implies  (h_1\otimes_1 h_2)\le_1 h_3\to_1 \psi(f) \\
&\implies \psi\varphi(h_1\odot_1 h_2)\le_1 h_3\to_1 \psi(f) \\
&\implies h_1\odot_1 h_2\le_1 h_3\to_1 \psi(f) \\
&\implies (h_1\odot_1 h_2)\odot_1 h_3 \le_1 \psi(f) \\
&\implies \psi\varphi((h_1\odot_1 h_2)\odot_1 h_3) \le_1 \psi(f) \\
&\implies \psi\varphi(h_1\odot_1 (h_2\odot_1 h_3)) \le_1 \psi(f) \\
&\implies h_1\otimes_1 (h_2\odot_1 h_3)) \le_1 \psi(f) \\
&\implies h_2\odot_1 h_3 \le_1 h_1\rightrightarrows_1 \psi(f) \\
&\implies \psi\varphi(h_2\odot_1 h_3) \le_1 h_1\rightrightarrows_1 \psi(f) \\
&\implies h_2\otimes_1 h_3 \le_1 h_1\rightrightarrows_1 \psi(f) \\
&\implies h_1\otimes_1 (h_2\otimes_1 h_3) \le_1  \psi(f)\\
&\implies h_1\odot_1 (h_2\otimes_1 h_3) \le_1  \psi(f).
\end{align*}
Thus $\psi(f)\in H_1 \implies \psi(f)\in H_2$.

\medskip
Conversely, let $f\in A_2^M$ such that $h_1\odot_1 (h_2\otimes_1 h_3)\leq_1\psi(f)$.
\begin{align*}
h_1\odot_1 (h_2\otimes_1 h_3)\leq_1 \psi(f) &\implies  h_2\otimes_1 h_3 \le_1 h_1\to_1 \psi(f) \\
&\implies \psi\varphi(h_2\odot_1 h_3)\le_1 h_1\to_1 \psi(f) \\
&\implies h_2\odot_1 h_3\le_1 h_1\to_1 \psi(f) \\
&\implies h_1 \odot_1 (h_2\odot_1 h_3) \le_1 \psi(f) \\
&\implies \psi\varphi(h_1\odot_1 (h_2\odot_1 h_3)) \le_1 \psi(f) \\
&\implies \psi\varphi((h_1\odot_1 h_2)\odot_1 h_3) \le_1 \psi(f) \\
&\implies (h_1\odot_1 h_2)\otimes_1 h_3 \le_1 \psi(f) \\
&\implies h_1\odot_1 h_2 \le_1 h_3\rightrightarrows_1 \psi(f) \\
&\implies \psi\varphi(h_1\odot_1 h_2) \le_1 h_3\rightrightarrows_1 \psi(f) \\
&\implies h_1\otimes_1 h_2 \le_1 h_3\rightrightarrows_1 \psi(f) \\
&\implies (h_1\otimes_1 h_2)\otimes_1 h_3 \le_1  \psi(f)\\
&\implies (h_1\otimes_1 h_2)\odot_1 h_3 \le_1  \psi(f).
\end{align*}
Thus $\psi(f)\in H_2 \implies \psi(f)\in H_1$. Therefore $H_1=H_2$, and
\[
h_1\otimes_1 (h_2\otimes_1 h_3) = (h_1\otimes_1 h_2)\otimes_1 h_3.
\]

\vspace*{-9mm}
\end{proof}

We can now define the product $\odot$ and the implication $\rightarrow$
on the set of all concepts.
\begin{theorem}\label{thm:main}
Let $\mathcal{A}_1$  and $\mathcal{A}_2$ be two complete residuated multilattices, $G$ and $M$ two non-empty sets
(of objects and attributes) and
$(\varphi, \psi)$ a Galois connection between $A_1^G$ and $A_2^M$. Then
\[\mathcal{C}=\{(h,f)\in A_1^G\times A_2^M; \varphi(h)=f \text{ and } \psi(f)=g \})\]
 is a complete residuated multilattice, if $\Ext(\mathcal{C})$ and
$\Int(\mathcal{C})$ are closed under $\to_1$ and $\to_2$ respectively, with
	\begin{align*}\label{e1}
	(h_1, f_1)\odot(h_2, f_2) &= (h_1\otimes_2 h_2, \varphi(h_1\otimes_2 h_2)) \text{ and } \\
	(h_1, f_1)\rightarrow(h_2, f_2) &= (h_1\rightrightarrows_2 h_2, \varphi(h_1\rightrightarrows_2 h_2)).
	\end{align*}
	for all concepts $(h_i, f_i)$, $i=1,2$.
\end{theorem}
\begin{proof}
	In (\cite{MOR13, MOPR16}) it is shown that $(\mathcal{C}, \preceq)$ is a complete multilattice.
	Let $(h_1, f_1)$, $(h_2, f_2)$ and $(h_3, f_3)$ in $\mathcal{C}$. From Lemma \ref{lemme 1} we know that
$(h_1\otimes_1
	h_2, \varphi(h_1\otimes_1 h_2))$ and $(h_1\rightrightarrows_2 h_2,
	\varphi(h_1\rightrightarrows_1 h_2))$ are concepts. It remains to
	prove that $\odot$ and $\rightarrow$ satisfy the adjointness
	condition that is,
	\begin{align*}
	(h_1, f_1)\odot(h_3, f_3)\le (h_2, f_2)&\iff (h_3, f_3)\le (h_1, f_1)\rightarrow(h_2, f_2).
\end{align*}
This is equivalent to
\begin{align*}
(h_1\otimes_1 h_3, \varphi(h_1\otimes_1 h_3))\le (h_2, f_2) &\iff
(h_3, f_3)\le (h_1\rightrightarrows_1 h_2, \varphi(h_1\rightrightarrows_1 h_2)),
\end{align*}
which is again equivalent to prove that
\[
h_1\otimes_1 h_3\leq_1
	h_2\Leftrightarrow h_3\leq_1 h_1\rightrightarrows_1 h_2.
\]
This is true since $\otimes_1$ and $\rightrightarrows_1$ satisfy the adjointness condition,  by
Lemma \ref{lemme 1}.
\end{proof}

We have thus proved that with residuated multilattices we can obtain a residuated concept  multilattice, assuming that the extents and intents are closed under $\to$. A special case is when the Galois connection $(\varphi,\varphi)$ preserves the residuation $(\odot,\to)$. By this we mean:
\begin{align*}
    \varphi(h_1\odot h_2) =  \varphi(h_1) \odot \varphi(h_2)\quad \text{ and }\quad \varphi(h_1\to h_2) =  \varphi(h_1) \to \varphi(h_2) \\
    \psi(f_1\odot f_2) =  \psi(f_1) \odot \psi(f_2)\quad \text{ and }\quad \psi(f_1\to f_2) =  \psi(f_1) \to \psi(f_2)
\end{align*}
for all $h_1,h_2\in A_1^G$ and $f_1,f_2\in A_2^M$.

\medskip
Replacing one of the residuated multilattices $\mathcal{A}_1$ and $\mathcal{A}_2$ by a residuated lattice, we
claim that the set of all concepts is a residuated lattice in the following corollary.
\begin{corollary} Under the assumption of the hypothesis of Theorem~\ref{thm:main}, if $\mathcal{A}_1$ or
$\mathcal{A}_2$ is a residuated lattice, then $(\mathcal{C}, \preceq)$ is a residuated lattice.
\end{corollary}
\begin{proof}
It was proved in \cite[Proposition 1]{RM12} that $(\mathcal{C}, \preceq)$ is a lattice.
Hence, by adding the above residuation, we have a residuated lattice.
\end{proof}

\section{Conclusion} \label{s:conclusion}
In this paper we have constructed a smallest (in term of cardinality) residuated multilattice that is not a residuated lattice. Using ordinal sums, we have shown how to produce residuated multilattices from old ones. Choosing residuated multilattices as set of truth values to evaluate objects and attributes, we have proved that the set of concepts forms a complete residuated multilattice whenever the Galois connection defining the concepts preserves the residuation. 

In \cite{MOR07} the authors discuss the use of ordered multilattices as underlying sets of truth-values for a generalized framework of logic programming. We believe that these results can be carried out to residuated multilattices, and plan to investigate these in our future work.

\end{document}